%% file: y-dec1123.tex
\documentclass{article}
\usepackage{latexsym}
\usepackage{amsmath}
\usepackage{amsfonts}
\usepackage{amssymb}
\input{Makro98}
\input{psfig}

\newcommand{\nn}{\nonumber \\}
\newcommand{\Idm}{\mbox{1\kern -0.5ex I}}

\newcommand{\openZ}{\mathbb{Z}}
\newcommand{\openC}{\mathbb{C}}
\newcommand{\openR}{\mathbb{R}}

\newcommand{\diag}{\mbox{\rm diag}}
\newcommand{\md}{\mbox{\rm mod }}

\newcommand{\bff}{{\bf f}}
\newcommand{\Id}{{\bf 1}}

\newcommand{\JJB}{\mathop{\JJ}_B} 
\newcommand{\JJA}{\mathop{\JJ}_A} 
\newcommand{\JJg}{\mathop{\JJ}_g} 

\newcommand{\JJBB}{\mathop{\JJ \,}\displaylimits_B} 
\newcommand{\JJbfBB}{\mathop{\JJ \,}\displaylimits_\bB} 
\newcommand{\JJAA}{\mathop{\JJ \,}\displaylimits_A} 
\newcommand{\JJgg}{\mathop{\JJ \,}\displaylimits_g} 

\newcommand{\Bprod}[2]{\mathop{#1 #2}\limits_B} 
\newcommand{\BBprod}[2]{\mathop{#1 #2}\displaylimits_B} 

\newcommand{\ed}{\end{document}} 
\begin{document}
\title{On the Decomposition of Clifford Algebras of Arbitrary Bilinear
Form\thanks{Paper presented at the 5th International Conference on
Clifford Algebras and their Applications in Mathematical Physics, Ixtapa,
Mexico, June 27 - July 4, 1999.}}
\author{Bertfried Fauser\\
Universit\"at Konstanz\\
Fakult\"at f\"ur Physik, Fach M678\\
D-78457 Konstanz\\
E-mail: Bertfried.Fauser@uni-konstanz.de
\and
Rafa\l \ Ab\l amowicz\\
Department of Mathematics Box 5054\\
Tennessee Technological University\\
Cookeville, TN 38505, USA\\
E-mail: rablamowicz@tntech.edu 
}
\date{September, 27 1999}
\maketitle

\begin{abstract}
Clifford algebras are naturally associated with quadratic forms. 
These algebras are $\openZ_2$-graded by construction. However, only 
a $\openZ_n$-gradation induced by a choice of a basis, or even better, 
by a Chevalley vector space isomorphism $\cl(V) \leftrightarrow 
\bigw V$ and an ordering, guarantees a multi-vector decomposition into
scalars, vectors, tensors, and so on, mandatory in physics. We show that
the Chevalley isomorphism theorem cannot be generalized to algebras 
if the $\openZ_n$-grading or other structures are added, e.g., 
a linear form. We work with pairs consisting of a Clifford algebra and 
a linear form or a $\openZ_n$-grading which we now call  
{\it Clifford algebras of multi-vectors\/} or 
{\it quantum Clifford algebras\/}. It turns out, that in this
sense, all multi-vector Clifford algebras of the same quadratic but
different bilinear forms are non-isomorphic. The usefulness of such
algebras in quantum field theory and superconductivity was shown
elsewhere. Allowing for arbitrary bilinear forms however spoils their 
diagonalizability which has a considerable effect on the tensor 
decomposition of the Clifford algebras governed by the periodicity 
theorems, including the Atiyah-Bott-Shapiro $\mod 8$ periodicity. We
consider real algebras $\cl_{p,q}$ which can be  decomposed in the
symmetric case into a tensor product $\cl_{p-1,q-1} \otimes \cl_{1,1}.$
The general case used in quantum field theory lacks this feature.
Theories with non-symmetric bilinear forms are however needed in the
analysis of multi-particle states in interacting theories. A connection to
$q$-deformed structures through nontrivial vacuum states in quantum
theories is outlined.

\noindent{\bf MSCS: 15A66; 17B37; 81R25; 81R50}

\noindent{\bf Keywords:} 
Clifford algebras of multi-vectors, Clifford map, quantum Clifford
algebras, periodicity theorems, index theorems, spinors, spin-tensors,
Chevalley map, quadratic forms, bilinear forms, deformed tensor products,
multi-particle geometric algebra, multi-particle states, compositeness,
inequivalent vacua
\end{abstract}

\section{Why study Clifford algebras of an arbitrary bilinear form?}

\subsection{Notation, basics and naming}

\subsubsection{Notation}

To fix our notation, we want to give some preliminary material. If
nothing is said about the {\it ring\/} linear spaces or algebras are
build over, we denote it by $\bR$ and assume usually that it is
unital, commutative and not of characteristic $2.$ In some cases we
specialize our base ring to the field of real or complex numbers denoted
as $\openR$ and $\openC.$

A {\it quadratic form\/} is a map $Q:V\mapsto \bR$ with the following
properties $(\alpha \in \bR,$ $\bV\in V)$
\begin{eqnarray}
i) && Q(\alpha \bV) \,=\, \alpha^2 Q(\bV), \nn
ii)&& 2g(\bx,\by) \,=\, Q(\bx-\by) - 
Q(\bx)- Q(\by),
\end{eqnarray}
where $g(\bx,\by)$ is bilinear and necessarily symmetric. $g({\bf
x},\by)$ is called {\it polar bilinear form\/} of $Q.$ Transposition
is defined as $g(\bx,\by)^T = g(\by,\bx).$ Quadratic forms
over the reals can always be diagonalized by a {\it choice\/} of a basis.
That is, in every equivalence class of a representation there is a 
diagonal representative.

We consider a {\it quadratic space\/} $\cH = (V,Q)$ as a pair of
a linear space $V$ --over the ring $\bR$-- and a quadratic form $Q.$
This is extended to a {\it reflexive space
\/} 
$\cH^\prime = (V,B)$ viewed as a pair of a linear space $V$ and an arbitrary
non-degenerate bilinear form $B = g+A,$ where $g=g^T$ and $A=-A^T$ are
the symmetric and antisymmetric parts respectively. $g$ is connected to a
certain~$Q.$

We denote the finite additive group of $n$ elements under addition modulo $n$
as $\openZ_n.$ This should not be confused with the ring $\openZ_n$ also 
denoted the same way.

Algebras or modules can be graded by an Abelian group. If the linear
space $W$ --not the same as $V$--, of an algebra can
be divided into a direct sum $W = W_0+W_1+\ldots +W_{n-1}$
and if the algebra product maps these spaces in a compatible way one onto 
another, see examples, so that the index labels behave like an Abelian
group, one refers to a {\it grading\/} \cite{Bourbaki}.\\

{\bf Example 1:} $W= W_0+W_1$ and $W_0W_0 \subseteq W_0,$ $W_0W_1 \simeq W_1W_0
\subseteq W_1$ and $W_1W_1 \subseteq W_0.$ The indices are added modulo 2 and
form a group $\openZ_2.$ If $W=W_0+W_1+\ldots+W_{n-1}$ one has e.g.
$W_iW_j \subseteq W_{i+j\,\, \mbox{\footnotesize \md}\, n}$ which is a
$\openZ_n$-grading. \\

In the case of $\openZ_n$-grading, elements of $W_m$ are
called $m$-vectors or homogenous multi-vectors. The elements of $W_0
\simeq \bR$ are also called {\it scalars\/} and the elements of $W_1$
are {\it vectors\/}. When the $\openZ_2$-grading is considered, one speaks
about even and odd elements collected in $W_0$ and $W_1$ respectively.

However, observe that the Clifford product is {\it not\/} graded in this
way since with $V \simeq W_1$ and $\bR \simeq W_0$ one has $V\times V
= \bR + W_2$ which is not group-like. Only the even/odd grading, 
sometimes called parity grading, is preserved, $\cl_+\cl_+ \subseteq
\cl_+,$ $\cl_+\cl_- \simeq \cl_-\cl_+ \subseteq \cl_-$ and $\cl_-\cl_-
\subseteq \cl_+.$ Hence, $\cl$ is $\openZ_2$ graded and $\cl \simeq \cl_+ + \cl_-
\simeq W_0 + W_1.$

Clifford algebras are displayed as follows: $\cl(B,V)$ is a 
{\it quantum Clifford algebra\/}, $\cl(Q,V)$ is a basis-free Clifford 
algebra, $\cl(g,V)$ is a Clifford algebra with a choice of a basis, 
$\cl_{p,q}$ and $\cl_n$ are real and complex Clifford algebras of 
symmetric bilinear forms with signature $p,q$ or of complex dimension 
$n$ respectively.

\subsubsection{Basic constructions of Clifford algebras}

Constructions of Clifford algebras can be found at various places in
literature. We give only notation and refer the Reader to these 
publications 
\cite{BennTucker,Bourbaki,BudinichTrautman,Chevalley,Crumeyrolle,Greub,Lounesto,Porteous}.\\

\noindent
{\bf Functorial:} 
The main advantage of the tensor algebra method is its
formal strength. Existence and uniqueness theorems are most easily
obtained in this language. Mathematicians derive almost all algebras from
the tensor algebra --the real mother of algebras-- by a process called
factorization. If one singles out a two-sided ideal $\cI$ of the tensor
algebra, one can calculate 'modulo' this ideal. That is all elements in
the ideal are collected to form a class called 'zero' $[0] \simeq \cI.$
Every element is contained in an equivalence class due to this construction.  
Denote the tensor algebra as $T(V)= \bR \oplus V \oplus \ldots
\otimes_n V\oplus \ldots$ and let $\bx,\by,\ldots \in V$ and $L,M,\ldots \in
T(V).$ This algebra is by construction naturally $\openZ_\infty$-graded
for any dimension of~$V.$

In the case of Clifford algebras, one selects an ideal of the form
\begin{equation}\label{cl-ideal}
\cI_{\cl} = \{ X \mid X=L\otimes(\bx\otimes \bx - Q(\bx)\Id)\otimes M\}
\end{equation}
which implements essentially the 'square law' of Clifford algebras. Note,
that elements of different tensor grades --scalar and grade two-- are
identified. Hence this ideal is not grade-preserving and the factor
algebra --the Clifford algebra-- cannot be $\openZ_\infty$-graded
(finiteness of $\cl(V))$; and not even multi-vector or $\openZ_n$-graded with
$n = \dim \, V$ because all indices are now $\md 2.$ However, the ideal
$\cI_{\cl}$ is $\openZ_2$-graded, that is, it preserves the evenness and the oddness
of the tensor elements. One defines now the Clifford algebra as:
\begin{equation}
\cl(Q,V) := \frac{\displaystyle T(V)}{\displaystyle \cI_{\cl}}.
\end{equation}

It is clear from the construction that a Clifford algebra is unital and
associative, a heritage from the tensor algebra.\\

\noindent
{\bf Generators and relations:}
Physicists and most people working in Clifford analysis prefer another
construction of Clifford algebras by generators and relations
\cite{CoxeterMoser}. One chooses a set of {\it generators} $\be_i,$ images
of some arbitrary basis elements ${\bx_i}$ of $V$ under the usual
Clifford map $\gamma : V \mapsto \cl(V)$ in the Clifford algebra
$\cl(V),$ and asserts the validity, in the case of $\bR = \openR$ or $\openC,$ of the normalized, 'square law':
\begin{equation}
\be_i^2 = \pm \Id.
\end{equation}
Using the linearity, that is polarizing this equations by $\be_i \mapsto
\be_i+\be_j,$ one obtains the usual set of relations which have
to be used to 'canonify' the algebraic expressions:
\begin{equation}
\label{eq:mod}
\be_i\be_j+\be_j\be_i = 2g(\be_i,\be_j)\Id = 2g_{ij}\Id.
\end{equation}
The definition of the Clifford algebra reads:
\begin{equation}
\cl(g_{ij},V) \simeq \mbox{\bf Alg}(\be_i)\, \md  \be_i\be_j \,=\, 2g_{ij}\Id -\be_j\be_i \,.
\end{equation}
While the --image of the-- numbers of the base field are called {\it scalars\/},
the $\be_i$ and their linear combinations are called {\it vectors\/}. The
entire algebra is constructed by multiplying and linear-combining the
generators $\be_i$ modulo the relation (\ref{eq:mod}). This 
'modulo relation' is in
fact nothing else as a 'cancellation law' which provides one with a {\it
unique\/} representative of the class of tensor elements. A basis of the
linear space underlying the Clifford algebra is given by {\it reduced
monomials\/} in the generators, where a certain ordering has to be chosen
in the index set, e.g. ascending indices or antisymmetry. A monomial
build out of $n$ generators and the linear span of such monomials is
called a homogenous $n$-vector. Thereby a unique $\openZ_n$-grading is
introduced
by the choice of a basis and an ordering.

This method has the advantage of being plain in construction, easy to
remember, and powerful in computational means.

\subsubsection{Naming}

A very important and delicate point in mathematics and physics is the
appropriate naming of objects and structures. Since we deal with a very
well known structure, but want to highlight special novel features, we have
to give distinguishing names to different albeit well known objects, which
otherwise could not be properly addressed. This section shall establish
such a coherent naming, at least for this article.

{\it Clifford algebra\/} is often denoted, following Clifford himself and
Hestenes, as 'Geometric Algebra', GA or 'Clifford Geometric Algebra' CGA or
'Clifford Grassmann Geometric Algebra' CGGA \cite{Parra-zac}. Having the
advantage of being descriptive this notation has, however, also a
peculiar tendency to call upon connotations and intuitions which
might {\it not in all cases\/} be appropriate. Even at this stage, one has
to distinguish 'Metric Geometric Algebra' MGA and 'Projective Geometric
Algebra' PGA which relies on the identification of the
homogenous multi-vector objects and geometrical entities
\cite{HestenesZiegler}. In the former case, 'vectors' are identified with
'places' --position vectors-- of pseudo-Euclidean or unitary spaces while
in the second case 'vectors' are identified with 'points' of a projective
space. 

{\it Both\/} variants, metric or projective, use unquestionably the
artificial multi-vector structure introduced by the mere notation of a
basis and foreign to Clifford algebras to assert 'ontological' statements
such as: '$\bx$ is a place in Euclidean space' or '$\bx$ is a
point in a projective space'.

{\it Both\/} of these interpretations have one thing in common, namely,
they assert an {\it object character\/} to the Clifford elements themselves. 
We will coin for this case the term '{\bf Classical Clifford Algebra}'.

To our current experience, the Wick isomorphism developed below
guarantees that such interpretation of Clifford algebras is {\it
independent\/} from the chosen $\openZ_n$-grading. That is, we make the
following {\it conjecture:\/} if the Clifford elements themselves are 
'ontologically' interpreted as 'place' or 'point' then all 
$\openZ_n$-gradings are isomorphic through the Wick isomorphism.

We turn to the second aspect. In \cite{Oziewicz97} Oziewicz introduced
the term 'Clifford algebras of multi-vectors' to highlight the fact that
he considered different $\openZ_n$-gradings or, equivalently, different
multi-vector structures. However, Clifford algebras have in nearly every
case been used as multi-vector Clifford algebras since mathematicians and
physicists want to consider the $n$-vectors or multi-vectors for different
purposes.

Following the introduction of Clifford algebras
of arbitrary bilinear forms, implicitly in \cite{Chevalley} and explicitly
in 
\cite{AblamowiczLounesto,Fauser-hecke,Fauser-thesis,Fauser-positron,Fauser-vacua,Fauser-q-groups,Fauser-mandel,Fauser-transition,Fauser-vertex,LounestoRiesz,Oziewicz-FGTC},
situations have occurred for good physical reasons where
different $\openZ_n$-gradings have led to different physical outcomes.
In those situations a theory of gradings is mandatory.

A new point is the {\it operational approach\/} to Clifford elements. If
one considers a Clifford number to be an operator, it has to act on
another object, a 'state vector'. This '{\it quantum point of view\/}'
moves also the ontological assertions into the states. Their
interpretation however is difficult.

Moreover, one has to deal with representation theory which was not
necessary in the 'classical' Clifford algebraic approach --in both senses
of classical, i.e. also as opposed to quantum, here. Adopting Wigner's
definition of a particle as an irreducible representation --of the
Poincar\'{e} group-- one has to seek irreducible representations of
Clifford algebras. It is a well known fact that these representations are
faithfully realized in spinor spaces. It is exactly at this place where
it will be shown in this article that one obtains different
$\openZ_n$-gradings or different multi-vector structures leading  to
different results. In fact we are able to find {\it irreducible\/} 
spinor spaces of dimension $8$ in $\cl_{2,2}(B,V),$ where $2,2$ denotes the 
signature of the symmetric part $g$ of $B,$ and not of dimension $4$ as 
predicted by the 'classical' Clifford algebra theory.

For the case of Clifford algebras of multi-vectors we coin the term {\bf
quantum Clifford algebra}.\footnote{This is close to Saller's notion of a 
``quantum algebra'' which denotes however a special choice of grading 
\cite{Saller-q-alg}.}$\phantom{}^{,}$\footnote{Classical Clifford algebras emerge as a
particular case of quantum Clifford algebras.} It is clear to us that we
risk creating a confusion with this term, which looks like a
$q$-deformed version of an ordinary Clifford algebra, while also in our case
the common 'square law' is fully valid! However, this link is not wrong!
As we show elsewhere in these proceedings \cite{AblaFauser}, one is able
to find Hecke algebras and $q$-symmetry {\it within\/} the structure of
the quantum Clifford algebra. It is also in accord with the attempt of G.
Fiore, presented at this conference, to describe $q$-deformed algebras in
terms of undeformed generators. This is just a reverse of our argument.
However, the characteristic point in our consideration is that we dismiss
the classical ontological interpretation in favor of an operational
interpretation. Thereby it is necessary to study states which are now
$\openZ_n$-grade dependent. Our approach should be contrasted by the
recent developments excellently described in \cite{Connes,Majid}. A
different treatment of Clifford algebras in connection with Hecke algebras
was given in \cite{Oziewicz95}.

As a last point, we emphasize that indecomposable spinor representations
of unconventionally large dimensions are expected to be spinors of bound
systems, see \cite{Fauser-vacua}. Hence, studying decomposability is the
first step towards an algebraic theory of compositeness including
stability of bound states.

\subsection{Why study $\cl(B,V)$ and not $\cl(Q,V)$? -- Physics}

Clifford algebras play without any doubt a predominant role in physics
and mathematics. This fact was clearly addressed and put forward by D.
Hestenes
\cite{Hestenes-spacetime,HestenesSobczyk,Hestenes-unifiedlang,Hestenes-foundmech}.
Based on this solid ground, we give an analysis of Clifford algebras of an
arbitrary bilinear form which exhibit novel features especially
regarding their representation theory. The most distinguishing fact between our
approach and usual treatments of Clifford algebras e.g., 
\cite{BennTucker,BudinichTrautman,Crumeyrolle,Lounesto,Porteous}, 
is that we seriously consider how the $\openZ_n$-grading is introduced in
Clifford algebras. This is most important since Clifford algebras are {\it
only\/} $\openZ_2$-graded by their natural --functorial-- construction.
The introduction of a further finer grading does therefore put new
assumptions into the theory. One might therefore ask, if theses
additional structures are important or even necessary in physics and
mathematics. 

Indeed, after examining various cases we notice that {\it every\/}
application of Clifford algebras which is computational --not only
functorial-- deals in fact with the so called Clifford algebras of
multi-vectors \cite{Oziewicz97} or {\it quantum Clifford algebras\/}.
However, the additional $\openZ_n$-grading, even if mathematically and
physically necessary for applications, is usually introduced without any
ado. Looking at literature we can however find lots of places where
$\openZ_n$-graded Clifford algebras are not only appropriate but needed.
This is in general evident in every quantum mechanical setup. 

If one analyzes functional hierarchy equations of quantum field theory
(QFT), one is able to translate these functionals with a help of Clifford
algebras. Such attempts have already been made by Caianiello
\cite{Caianiello}. He noticed that at least two types of orderings are
needed in QFT, namely the time-ordering and normal-ordering. Since one has
--at least-- two possibilities to decompose Clifford algebras into basis
monomials, he introduces Clifford and Grassmann bases. A basis of a
Clifford algebra is usually given by monomials with totally ordered index
sets. If one has a finite number of 'vector' elements $\be_i,$ one can, by
using the anti-commutation relations of the Clifford algebra, introduce the
following bases
\begin{eqnarray}
i) && \{\Id;\be_1,\ldots,\be_n;\be_1\be_2,\ldots;
\be_{i_1}\be_{i_2}\be_{i_3}{}_{(i_1<i_2<i_3)},\ldots\} \nn
ii) && \{\Id;\be_1,\ldots,\be_n;\be_{[1}\be_{2]},\ldots;
\be_{[i_1}\be_{i_2}\be_{i_3]},\ldots\} . 
\end{eqnarray}
We used the $[\ldots]$ bracket to indicate antisymmetrization in the index
set. An ordering of index sets is inevitable since the $\be_i\be_j$ and
$\be_j\be_i$ monomials are not algebraically independent due to the
anti-commutation relations $\be_i\be_j=-\be_j\be_i+g_{ij}\Id.$ Caianiello
identifies then the two above choices with time- and normal-ordering.
However, already at this point it is questionable why one uses
'lexicographical' ordering '$<$' and not e.g. the 'anti-lexicographical'
ordering '$>$' or an ordering which results from a permutation of the
index set.

A detailed study shows that fermionic QFT needs antisymmetric
index sets and that there are infinitely many such choices
\cite{Fauser-thesis,Fauser-transition}. Using this fact we have been able
to show that singularities, which arise usually due to the reordering
procedures such as the normal-ordering, are no longer present in such algebras 
\cite{Fauser-vertex}. Studying the transition from operator
dynamics to functional hierarchies, the so-called Schwinger-Dyson-Freese
hierarchies, in \cite{Fauser-thesis,Fauser-transition} it turned out that the
multi-vector structure, or, equivalently a uniquely chosen $\openZ_n$-grading,
was a {\it necessary input\/} to QFT.

Multi-particle systems provide a further place where a careful study of
gradings will be of great importance. It is a well known fact that one
has the Clebsch-Gordan decomposition of two spin-$\frac12$ particles
as follows \cite{FultonHarris,Hamermesh}:
\begin{equation}
\frac{{\bf 1}}{{\bf 2}}\otimes\frac{{\bf 1}}{{\bf 2}}
=
{\bf 0} \oplus {\bf 1}\,.
\end{equation} 
However, since this is an identity, it can be used either from left to right
to form bosonic spin $0$ and spin $1$ 'composites' {\it or\/} from right to
left! There is no way --besides the experience-- to distinguish if such a
system is composed, that is, dynamically stable or not, see
\cite{Fauser-positron}. From a mathematical point of view one cannot
distinguish $n$ free particles from an $n$-particle bound system by means
of algebraic considerations. This is seen clearly in the decomposition
theorems for Clifford algebras where larger Clifford algebras are
decomposed into smaller blocks of Clifford tensor factors. This cannot be
true for bound objects which lose their physical character when being
decomposed. An electron and proton system is quite different from a
hydrogen atom. In this work, we will see, that one can indeed find such
{\it indecomposable\/} states in quantum Clifford algebras.

This raises a question how to distinguish such situations. One knows from
QFT that interacting systems have to be described in non-Fock states and
that there are infinitely many such representations \cite{Haag-theorem}. It
is thus necessary to introduce the concept of {\it inequivalent states\/}
in finite dimensional systems \cite{Fauser-vacua,Ker1,Ker2}. Such states
are necessarily non-Fock states, since Fock states belong to
systems of non-interacting particles. This is the so-called {\it free case\/}
which is however very useful in perturbation theory. The present paper
supports the situation found in \cite{Fauser-vacua}.

Closely related to these inequivalent states are {\it condensation
phenomena\/}. As it was shown in \cite{Fauser-vacua}, one can algebraically
determine boundedness using an appropriate $\openZ_n$-grading. Furthermore, it
was shown that the dynamics {\it determines\/} correct grading. In BCS
theory of superconductivity the fact that bound states can or cannot be
build was shown to imply a gap-equation \cite{Fauser-vacua} which
governs the phase transition.

A further point related to $\openZ_n$-graded Clifford algebras is
$q$-quantization. This can be seen when studying physical systems as in
\cite{Fauser-fkt} and when adopting a more mathematical point of view 
as in \cite{Fauser-hecke,Fauser-q-groups}. In these proceedings a detailed
example was worked out to show how $q$-symmetry and Hecke algebras can be
described within quantum Clifford algebras \cite{AblaFauser}. It is quite
clear that this structure should play a major role in the discussion of
the Yang-Baxter equation, the knot theory, the link invariants and in other 
related fields which are crucial for the physics of integrable systems in 
statistical physics.

However, the most important implication from these various applications
is that the $q$-symmetry and more general deformations are {\it symmetries of
composites\/}. This was already addressed in \cite{Fauser-fkt} and more
recently in \cite{Fauser-hecke}. Also the present work provides full
support for this interpretation, as the talk of G. Fiore at this
conference. Providing as much evidence as possible to this fact was a
major motivation for the present work. 

\subsection{Why study $\cl(B,V)$ and not $\cl(Q,V)$? -- Mathematics}

There are also arguments of purely mathematical character which force us
to consider quantum Clifford algebras.

If we look at the construction of Clifford algebras by means of the tensor
algebra, we notice that $\cl$ is a functor. To every quadratic space
$\cH = (V,Q),$ a pair of a linear space $V$ over a ring $\bR$ and a
quadratic form $Q,$ there is a uniquely connected Clifford algebra
$\cl(Q,V).$ That is, one can introduce the algebra structure without any
further input or choices, so to say for free. One may further note that
if the characteristic of the ring $\bR$ is not $2,$ 
then there is a one-to-one correspondence between
quadratic forms and classes of symmetric matrices \cite{SchejaStorch}.
In other words, every symmetric matrix is a representation of a quadratic
form in a special basis. Over the reals (complex numbers) the classes of
quadratic forms can be labeled by dimension $n$ and signature $s$
(dimension $n$ only, no signature in $\openC).$ Equivalently one can use
the numbers $p,q$ of positive and negative eigenvalues of the quadratic
form. This leads to a classification (naming) of real (and complex)
Clifford algebras. One writes $\cl(Q,V) \simeq \cl_{p,q}$ $(\cl_n)$
where $\dim \,V=n=p+q$ and $Q$ has signature $s=p-q.$ The remarkable fact is that
the 'square law' for vectors $Q(\bv) \equiv \bv^2 = \alpha
\Id \in \cl(Q;V)$ $(\alpha \in \openR$ or $\alpha \in \openC)$
is a diagonal map determining only the symmetric part of the map
$Q(\bV)\mapsto \openR.$ Following Clifford one should note that
the product operation can be seen as 
acting on the second factor $2 \times x$ as a doubling of $x;$ that is, $2 \, \times $ 
is a doubling operator or endomorphism acting on the space of the 
second factor. In this sense any 'Clifford number' induces an
endomorphism on the graded space $W$ underlying the algebra and it is
questionable why one should use only diagonal maps and their symmetric
polarizations. Furthermore, note that one has
\begin{eqnarray}
\mbox{quadratic forms} &\simeq& \frac{\mbox{bilinear
forms}}{\mbox{alternating forms}} .
\end{eqnarray} 
The dualization $V\mapsto V^* \simeq \mbox{lin-Hom}(V,\openR)$ is
performed by an arbitrary (non-degenerate) bilinear form. Endomorphisms
have in general the following form
\begin{equation}
\mbox{End}(V) \simeq  V\otimes V^*,
\end{equation}
so why do we restrict ourselves to the symmetric case? If we consider a
pair $(V,B)$ of a space $V$ and an arbitrary bilinear form $B,$ can we
construct functorially an algebra like the Clifford algebra for the pair
$\cH = (V,Q)$?

It can be easily checked that if one insists on the validity of the
'square law' $\bv^2=\alpha\Id,$ the {\it anti-commutation
relations\/} of the resulting algebra are the same as for usual Clifford
algebras while the {\it commutation relations\/} --and thus the meaning of
ordering and grade-- is changed. Let $B=g+A,$ $A^T=-A,$ $g^T=g.$ We
denote $B(\bx,\by) = \bx \JJBB \by,$ $A(\bx,\by) = \bx \JJAA \by$ and 
$g(\bx,\by) = \bx \JJgg \by$ (the latter also denoted by Hestenes and
Sobczyk as $\bx \cdot \by).$\footnote{The symbols $\JJBB,\,\JJAA$ and $\JJgg$ denote the left contraction in $\cl(B,V)$ with respect to $B,\,A$ and $g$ respectively.} Then, the $B$-dependent Clifford product 
$\Bprod{\bx}{\by}$ of two $1$-vectors $\bx$ and $\by$ in $\cl(B,V)$ can be
decomposed in {\it different ways\/} into scalar and bi-vector parts as follows
\begin{eqnarray}\label{eq-11}
\BBprod{\bx}{\by} = \bx \JJg \by + \bx \dw \by  && \mbox{Hestenes, common case,
$A=0$} \nn
\BBprod{\bx}{\by} = \bx \JJBB \by + \bx \w  \by && \mbox{Oziewicz, Lounesto, Ab\l amowicz, Fauser},   
\end{eqnarray}
where $\bx \dw \by = \bx \w \by + A(\bx, \by) = \bx \w \by + \bx \JJA \by.$ Of course, 
for any $1$-vector $\bx$ and any element $u$ in $\cl(B,V)$ we have:
\begin{equation}\label{eq-12}
\BBprod{\bx}{u} = \bx \JJBB u + \bx \w u = 
\bx \JJgg u + \bx \JJAA u + \bx \w u = \bx \JJgg u + \bx \dw u.  
\end{equation}
Notice that the element $\bx \dw u = \bx \JJA u + \bx \w u$ is not even a 
homogenous multi-vector in $\bigw V.$ We have thus established that
the multi-vector structure is uniquely connected with the antisymmetric
part $A$ of the bilinear form, see also
\cite{AblamowiczLounesto,Fauser-mandel,Fauser-transition}.

This has an immediate consequence: in some cases one finds bi-vector
elements which satisfy minimal polynomial equations of the Hecke type
\cite{Fauser-hecke,Fauser-q-groups}. This feature is treated extensively
elsewhere in this Volume \cite{AblaFauser}.

Some mathematical formalisms, not treated here, are closely connected to this
structure. One is the structure theory of Clifford algebras over arbitrary
rings \cite{Hahn} where a classification is still lacking. Connected to
these questions is the arithmetic theory of Arf invariants and the
Brauer-Wall groups.

Much more surprising is the fact that due to central extensions the
ungraded bi-vector Lie algebras turn into Kac-Moody and Virasoro algebras
\cite{MoodPianzola} and, as it is also shown in \cite{AblaFauser}, to
some $q$-deformed algebras.

Since Clifford algebras naturally contain reflections, automorphisms
generated by non-isotropic vectors, we expect to find infinite
dimensional Coxeter groups \cite{CoxeterMoser,Hiller}, affine Weyl groups
etc., connected to $\openZ_n$-graded or quantum Clifford algebras.

Involutions connected to special elements, norms and traces \cite{Hahn} 
are also affected by different gradings. This has considerable effects.
One important point is that the Cauchy-Riemann differential equations are
altered which makes probably the concept of monogeneity
\cite{Lounesto} grade dependent. However, this is speculative.

\section{Chevalley's approach to Clifford algebras}

\subsection{Confusion with Chevalley's approach}

Chevalley's book {\it ``The algebraic theory of spinors''\/}
\cite{Chevalley} seems to have been badly accepted by working mathematicians
and physicists despite its frequent citation. Albert Crumeyrolle stated the 
following in \cite{Crumeyrolle}, p. xi:
\begin{quote}
{\it In spite of its depth and rigor, Chevalley's book
proved too abstract for most physicists and the notions
explained in it have not been applied much until recently,
which is a pity.\/}
\end{quote}
The more compact and readable book {\it ``The study of certain important
algebras''\/} \cite{Chevalley2} seems to be little known. However, one can
find in many physical writings e.g. Berezin \cite{Berezin} very analogous
structures, without mentioning the much more complete work of
Chevalley.

When looking for the most general construction of Clifford algebras over
arbitrary rings including the case where the characteristic of $\bR$
is $2,$ Chevalley constructed the so-called {\it Clifford map.\/} This map is an
injection of the linear space $V$ into the algebra $\cl(V)$ which establishes
the 'square law'. This construction emphasizes the operator character of
Clifford algebras and establishes a connection between the
spaces underlying the $\openZ_n$-graded Grassmann algebra and the thereon
constructed Clifford algebra. For our purpose it is important
that {\it only\/} Chevalley's construction allows a non-symmetric bilinear
form in constructing Clifford algebras. However, this fact is not
explicit in Chevalley's writings but it is clearly emphasized in
\cite{Oziewicz-FGTC}.

Ironically, a careful analysis of Lounesto shows that even Crumeyrolle
made a mistake in describing the Chevalley isomorphism connecting
Grassmann and Clifford algebra spaces. In \cite{Lounesto-counter} Lounesto
points out that Crumeyrolle rejects the Chevalley isomorphism for {\it
any\/} characteristic. This seems to be implied by Crumeyrolle's
frequent questioning, also in previous Clifford conferences of this
series: ``What is a bi-vector?'' \cite{Crumeyrolle-AL}. However, an
isomorphism can be uniquely given if the characteristic of $\bR$ is
not $2,$ see \cite{LounestoRiesz,Lounesto-counter}. On the other hand,
Lounesto points out that Lawson and Michelsohn \cite{LawsonMichelsohn}
postulate such an isomorphism which is wrong in the exceptional case of
characteristic $2.$ One should note in this context that their point of
view is taken by almost all working mathematicians and physicists.

At this point we submit, that we insist on Chevalley's construction even
in the case of characteristic not $2.$ Lounesto claims that in this cases
$\cl(B,V)$ is isomorphic to $\cl(Q,V)$ with $Q$ the quadratic form
associated to $B.$ In fact, this is true for the Clifford algebraic
structure and was proved in \cite{AblamowiczLounesto} up to the dimension $9$
of $V.$ However, this, the so-called {\it Wick isomorphism\/} between $\cl(B,V)$
and $\cl(Q,V),$ has to be rejected when the $\openZ_n$-grading is
considered, or, in other words, the multi-vector structure. Hence, we
reject Lounesto's judgment that it is worth studying $\cl(B,V)$ only in
characteristic $2$ for the reason of carefully treating the involved
$\openZ_n$-grading or multi-vector structure. This is one of the main
points of our analysis.

\subsection{Chevalley's construction of $\cl(B,V)$}

A detailed and mathematical rigorous development of quantum Clifford
algebras $\cl(B,V)$ can be found in \cite{Fauser-hecke,Fauser-transition}.
We will develop only the notation and point out some peculiar features
insofar as they appear in the present study, see also 
\cite{AblamowiczLounesto,Fauser-mandel}.

The main feature of the Chevalley approach is that Clifford algebras are
constructed as special --satisfying the 'square law'-- endomorphism
algebras on --the linear space of-- a Grassmann algebra. In this way the
Grassmann algebra, which is naturally $\openZ_n$-graded, induces via the
Chevalley isomorphism a grading or multi-vector structure in the Clifford
algebra. This grading is however not preserved by the Clifford product
which renders the Clifford algebra to be a deformation of the Grassmann
algebra.

To proceed along this line we construct the Grassmann algebra as a factor
algebra of the tensor algebra. Let
\begin{equation}
\cI_G := \{ X \mid X=A\otimes(\bx\otimes \bx)\otimes B \}
\end{equation}
with notation as in (\ref{cl-ideal}) and define
\begin{equation}
\bigwedge V := \frac{T(V)}{\cI_G},
\quad
\pi : T(V) \,\mapsto\, \bigwedge V\,.
\end{equation}
The projected tensor product $\pi(\otimes) \mapsto \w$ is denoted as
wedge or outer product. The induced grading is
\begin{equation}
\bigwedge V = \bR \oplus V\w V \oplus \ldots \oplus \w^n V \oplus
\ldots\quad .
\end{equation}

As the next step, we consider reflexive duals of the linear space $V.$ Define
\begin{equation}
V^* := \mbox{lin-Hom}(V,\bR)
\end{equation}
where $\dim \, V^* = \dim \,V$ (reflexivity). Using the action of
the dual elements on $V$ we define the (left) contraction $\JJBB$ as:
\begin{equation}
i_\bx(\by) = \bx\JJBB \by \,=\, B(\bx,\by).
\end{equation}
Note, that $i_\bx\in V^*$ is the dualized element $\bx$ and that
here a {\it certain\/} duality map is employed. If this is the usual
duality map $i_{\be_i}(\be_j) = \delta_{ij}$ one denotes this as Euclidean
dual isomorphism and writes the map as $\star$ \cite{Saller-q-alg}. The
notation $\bx\JJB \by$ and much more $B(\bx,\by)$ is very
peculiar since we have
\begin{eqnarray}
\JJBB \,:\, V\times V \mapsto V, &&
B     \,:\, V\times V \mapsto V. 
\end{eqnarray} 
Hence, $\JJBB$ and $B$ are in $\mbox{lin-Hom}(V\times V, \bR) \simeq 
V^* \times V^*.$ In this notation a dual isomorphism is {\it
implicitly\/} involved, since we consider really maps of the form
\begin{equation}
<.\mid.> \,:\, V^* \times V \mapsto \bR
\end{equation}
which might be called {\it a dual product\/} or {\it a pairing\/} 
\cite{Bourbaki,Saller-q-alg}. 

Having defined the action of $V^*$ on $V,$ we lift this action to the
entire Grassmann algebras $\bigwedge V$ and $\bigwedge V^*.$ For 
$\bx,\by\in V,$ and $u,v,w\in \bigwedge V$ we have:
\begin{eqnarray}
i)  && \bx \JJBB \by = B(\bx,\by), \nn
ii) && \bx \JJBB (u\wedge v) = (\bx \JJBB u) \wedge v
       + \hat{u} \wedge (\bx \JJBB v), \nn
iii)&& (u \wedge v) \JJBB w = u \JJBB (v \JJBB w),
\end{eqnarray}
where $\hat{~}$ is the involutive map --grade involution-- $\hat{~} : V
\mapsto -V$ lifted to $\bigwedge V.$ The Clifford algebra $\cl(B,V)$ is
then constructed in the following way. Define an operator 
$L^\pm_\bx : \bigwedge V \mapsto \bigwedge V$ for any $\bx \in V$ as:
\begin{equation}
(L^\pm_\bx)^2 := \bx \JJBB \cdot \, \pm \, \bx \w \cdot
\end{equation}
and observe that this is a Clifford map
\cite{Chevalley,Fauser-hecke,Fauser-transition}
\begin{equation}
(L^\pm_\bx)^2 = \pm Q(\bx) \Id,
\end{equation}
where $Q(\bx) = B(\bx,\by).$ This is nothing else as again 
the 'square law', and one proceeds as in the case of generators and
relations. Chevalley has thus established that
\begin{eqnarray}
\cl(B,V) \subset \mbox{End}(\bigwedge V).
\end{eqnarray}
This inclusion is strict.

\section{Wick isomorphism and $\openZ_n$-grading}
\subsection{Wick isomorphism}

In this section we will prove the following {\bf Theorem:}
\begin{eqnarray}
\cl(B,V) \cong \cl(Q,V)
\end{eqnarray}
as $\openZ_2$-graded Clifford algebras.

This isomorphism, denoted below by $\phi,$ is the {\it Wick isomorphism\/}
since it is the well know normal-ordering transformation of the quantum field
theory \cite{Dyson,Fauser-transition,Wick}. This was not noticed for a
long time which is another 'missed opportunity' \cite{Dyson-mo}.\\

\noindent{\bf Proof:}
The proof proceeds in various steps, numbered by letters a, b, c, etc. After
defining the outer exponential, we prove the following important formulas:
\begin{eqnarray}\label{eqs-to-prove}
i)   && e_\wedge^{-F}\wedge e_\wedge^{F} = \Id, \nn
ii)  && e_\wedge^{-F}\wedge \bx \wedge e_\wedge^{F} \wedge u 
\,=\, \bx\wedge u, \nn
iii) && e_\wedge^{-F}\wedge (\bx \JJgg (e_\wedge^{F} \wedge u)) 
\,=\, \bx\JJgg u + (\bx \JJgg F) \wedge u,
\end{eqnarray}
and finally we show that the Wick isomorphism $\phi$ is given as:
\begin{eqnarray}
\cl(B,V)   &=& \phi^{-1} (\cl(g,V)) \nn
           &=& e_\wedge^{-F}\wedge \cl(Q,V)\wedge e_\wedge^F \\
       &\cong& (\,\cl(g,V),<.>^A_r \,)
\end{eqnarray}
where $<.>^A_r$ denotes the $A$-dependent $\openZ_n$-grading.

That is, the isomorphism is given by the following transformation
of {\it vector\/} variables which is then algebraically lifted to 
the entire algebra:
\begin{eqnarray}\label{subs-law}
\bx \JJgg \cdot   & \rightarrow& \bx \JJBB \cdot \,=\, \bx \JJgg \cdot \, + \,
                 (\bx \JJgg  F)\wedge \cdot \nn
\bx \wedge \cdot  & \rightarrow& \bx \wedge \cdot 
\end{eqnarray}
\noindent{\bf a)}
According to Hestenes and Sobczyk \cite{HestenesSobczyk} it is possible to
express every antisymmetric bilinear form in the following way
\begin{equation}
A(\bx,\by) := F\JJgg (\bx\w \by)
\end{equation}
where $F$ is an appropriately chosen bi-vector. $F$ can be decomposed in a
non-unique way into homogenous parts $F_i= \ba_i\w \bb_i,$ $F=\sum
F_i.$ We define the outer exponential of this bi-vector as $(\w^0 F = \Id)$
\begin{eqnarray}
e_\w^{F} := \sum \frac{1}{n!}\w^n\,F \,=\,
\Id+F+\frac{1}{2}F\w F+ \ldots + \frac{1}{n!}\w^n \, F +\ldots\,\,. 
\label{expans}
\end{eqnarray}
This series is finite when the dimension of $V$ is finite since in that
case there exists a term of the highest grade.

\noindent{\bf b)}
Substitute the series expansion (\ref{expans}) into (\ref{eqs-to-prove}-i) and note that after
applying the Cauchy product formula for sums we have
\begin{equation}
e_\wedge^{-F}\wedge e_\wedge^{F} = \sum_{r=0}^\infty \left(\sum_{l=0}^r (-1)^l \binom{r}{l}\right) \frac{1}{r!}\w^r F.
\end{equation}

The alternating sum of the binomial coefficient is zero except in the case
$r=0$ when we obtain \Id, which proves formula (\ref{eqs-to-prove}-i).

\noindent{\bf c)}
To prove (\ref{eqs-to-prove}-iii) one needs the commutativity of 
$\bx \JJgg F_i$ with $F_j.$ If the contraction is zero, it commutes
trivially, if not, the contraction is a vector \by. From $\by \w
F = F \w \by$ for every bi-vector, we have that $\bx\JJgg F_i$
commutes with any $F_j$ and thus with $F.$ This allows us to write
\begin{equation}
\bx \JJgg (\w^n F) = n (\bx \JJgg F) \w (\w^{(n-1)} F). 
\end{equation}
Once more using $\w^0 F = \Id,$ the Leibniz' rule and the fact
that $\hat{F}=F,$ we obtain
\begin{equation}
\bx \JJgg \,(e_\w^F \w u) = e_\w^F \w (\bx \JJgg u + (\bx \JJgg F) \w u ),
\end{equation}
which proves (\ref{eqs-to-prove}-iii).

\noindent{\bf d)}
Since any vector $\by$ commutes under the wedge with any bi-vector $F,$
the case (\ref{eqs-to-prove}-ii) reduces to {\bf b)}.

\noindent{\bf e)}
The Wick isomorphism is now given as $\cl(B,V)=\phi^{-1}(\cl(g,V))
= e_\wedge^{-F}\wedge \cl(Q,V)\wedge e_\wedge^F.$ The same transformation
can be achieved by decomposing every Clifford 'operator' into vectorial
parts and then into contraction and wedge parts w.r.t. $(g,\w)$ and then
performing the substitution laws given in (\ref{subs-law}) and a final
renaming of the contractions; see \cite{Fauser-transition} for
an application in quantum field theory. 

Note, that since the wedges are {\it not\/} altered and the
new contractions are given by $\bx \JJBB \cdot \equiv \mbox{d}_\bx (\, \cdot \, ) := 
\bx \JJgg \cdot \, + \, (\bx \JJgg F) \w \cdot,$ this transformation does mix
grades, but it respects the parity. It is thus a $\openZ_2$-graded
isomorphism. QED. \\

An equivalent proof was delivered in \cite{StumpfBorne} without using
(explicitly) Clifford algebras but index doubling --see below. The
Wick isomorphism was called there 'nonperturbative normal-ordering'.

\subsection{$\cl(B,V) \,\leftrightarrow\, \cl(Q,V)$ -- Isomorphic yet
different?}

We have already discussed that many researchers reject the idea that
$\cl(B,V)$ is of any use because of the Wick isomorphism. However, as our
proof has shown this isomorphism is only $\openZ_2$-graded. Indeed it was
not the mathematical opportunity, but a necessity in modeling quantum
physical multi-particle systems and quantum field theory which forced us
to investigate quantum Clifford algebras 
\cite{Fauser-thesis,Fauser-vacua,Fauser-transition,Fauser-vertex}.

Decomposing $B$ into $g,A$ as in (\ref{eq-12}) and noting that in our
case, of characteristic not $2$ one has $Q(\bx) = g(\bx,\bx),$ one
concludes that $\cl(Q,V)$ is exactly the equivalence class of
$\cl(B,V)\simeq \cl(g+A,V)$ with $A$ varying arbitrarily:
\begin{equation}
\cl(Q,V) = [ \cl(g+A,V) ] \, .
\end{equation}
In other words, one does not have a single Clifford algebra $\cl(Q,V)$ but
an entire class of equivalent --under the $\openZ_2$-graded
Wick isomorphism-- Clifford algebras $\cl(B,V).$ This can be written as
\begin{equation}
\cl(Q,V) \simeq \cl(g+A,V)\; \md A 
\end{equation}
which induces a unique projection from the class of quantum Clifford 
algebras onto the classical Clifford algebra. Such a projection $\pi$ can
be defined as:
\begin{eqnarray}
i)  && \pi \,:\, T(V) \mapsto \cl(B,V) \nn
ii) && <.>^A_r := \pi(\otimes^r\,V).
\end{eqnarray}
This is once more a sort of 'cancellation law'. The important fact is that
{\it only\/} those properties belong to $\cl(Q,V)$ which do {\it not\/}
depend on the particular choice of a representant parameterized by $A.$ 
Physically speaking, only those properties belong to $\cl(Q,V)$ which are 
homogenous over the entire equivalence class.

As we will show now, especially the multi-vector $\openZ_n$-grading is {\it
not\/} of this simple type. Recall that it is possible to decompose the 
Clifford product in various ways as in (\ref{eq-11}) and (\ref{eq-12}). 
Hence we obtain a relation between the $\w$- and the $\dw$-grading as:
\begin{equation}
\bx \dw \by = A(\bx,\by) + \bx \w \by
\end{equation}
which shows that a $\dw$-bi-vector is an inhomogeneous $\w$-multi-vector
and vice versa. Since the antisymmetric part can be absorbed in the wedge
product, using the Wick isomorphism, we can give the grading explicitly
by writing
\begin{equation}
<.>^A_r \;=\; <.>^{\dw}_r
\end{equation}
with respect to the doted wedge $\dw$ {\it within\/} the undeformed
algebra
$\cl(Q,V),$ see G. Fiore's talk. This gives us a second
characterization of $\cl(B,V),$ namely
\begin{equation}
\cl(B,V) \simeq (\cl(Q,V),<.>^A_r) \, .
\end{equation}
That is, $\cl(B,V)$ can be seen as a pair of a classical
$\openZ_2$-graded Clifford algebra $\cl(Q,V)$ and a unique multi-vector
structure given by the projectors $<.>^A_r.$ As a main result we have that
these algebras are {\it not\/} isomorphic under the Wick isomorphism
\begin{eqnarray}
\cl(g+A_1,V) &{\mathop{\not\simeq}\limits_{\mbox{\footnotesize{Wick}}}}& \cl(g+A_2,V) 
\quad \mbox{iff}\quad A_1 \,\not=\, A_2 \, .
\end{eqnarray}

\section{Periodicity theorems}

Our theory will have an impact on all famous periodicity theorems of
Clifford algebras, especially on the Atiyah-Bott-Shapiro $\md 8$
index theorem \cite{AtiyahBottSharpiro}. But to be as concrete and explicit as
possible, we restrict ourself to the case $\cl_{p,q} \simeq
\cl_{p-1,q-1} \otimes \cl_{1,1}.$ Periodicity theorems can be found, for example, in 
\cite{BennTucker,BudinichTrautman,Lam73,Maks,Porteous}.

We need some further notation. Let $V_{p,q} = (g_{p,q},V)$ be a quadratic
space, where $g = \diag(1,\ldots,1,-1,\ldots,-1)$ with $p$ plus signs and $q$ minus
signs, and let $V$ be a linear space of dimension $p+q.$ According to the
Witt theorem \cite{Witt-thrm} one can split off a quadratic space 
of the hyperbolic type $M_{1,1}.$ This split is orthogonal with respect to $g:$
\begin{equation}
V_{p,q} = N_{p-1,q-1} \perp_g M_{1,1} \, .
\end{equation}
If one applies the Clifford map $\gamma : V_{p,q} \mapsto \cl_{p,q}$ and
defines its natural restrictions $\gamma^\prime : N_{p-1,q-1} \mapsto
\cl_{p-1,q-1},$ $\gamma^{\prime\prime} : M_{1,1} \mapsto \cl_{1,1},$ one
obtains the following {\bf Periodicity Theorem:}
\begin{equation}
\cl_{p,q} \simeq \cl_{p-1,q-1} \otimes \cl_{1,1}.
\end{equation}
While in this special case the tensor product may be ungraded, in general the 
tensor product in such decompositions may be graded or not, see
\cite{BudinichTrautman,Lam73,Maks}.

Using the obvious notation $\cl(V_{p,q}) = \cl_{p,q}(Q)$ and
introducing the restrictions of the Wick isomorphism $\phi^{-1}\vert_N$ and
$\phi^{-1}\vert_M,$ (here $N=N_{p-1,q-1}$ and $M=M_{1,1}),$ we can calculate the 
decomposition of $\cl_{p,q}(B).$ However, if there are terms in the bi-vector $F$ 
which connect spaces $N$ and $M,$ that is, if $F=\sum F_i$ and if there exists  
$F_s=\ba_s \wedge {\bf b}_s$ with $\ba_s \in N,$ $\bb_s \in M,$ this part of the
construction belongs {\it neither\/} to the restriction $\phi^{-1}\vert_N$
{\it nor\/} to $\phi^{-1}\vert_M.$ We have {\it either\/} {\bf no tensor
decomposition} {\it or\/} a {\bf deformed tensor product}. Expressed in
formulas we get:

\begin{eqnarray}
\cl_{p,q} &=& \phi^{-1}(\cl_{p,q}(Q)) \nn
&=&
\phi^{-1}\left[
\cl_{p-1,q-1}(Q\vert_N)\otimes\cl_{1,1}(Q\vert_M)
\right] \nn
&=&\cl_{p-1,q-1}(B\vert_N)(\phi^{-1}\otimes)\cl_{1,1}(B\vert_M) \nn
&=&\cl_{p-1,q-1}(B\vert_N)\,\otimes_{\phi^{-1}}\, \cl_{1,1}(B\vert_M).
\end{eqnarray}

\noindent{\bf Remark:}
The deformed tensor product $\otimes_{\phi^{-1}}$ is not braided by
construction, since we have no restrictions on $\phi^{-1}.$ But one is 
able to find e.g., Hecke elements, etc., necessary
for a common $q$-deformation or, more generally, a braiding.

As the main result of our investigation we have shown that quantum
Clifford algebras {\it do not\/} come in general with periodicity theorems
as e.g. the famous Atiyah-Bott-Shapiro $\md 8$ index theorem. This has
{\it enormous\/} impact on quantum manifold theory and the topological
structure of such spaces as well as on their analytical properties. However, we
have constructed a deformed --not necessarily braided-- tensor product
$\otimes_{\phi^{-1}}$ which gives a decomposition at the cost of losing
(anti)-commutativity. To fully support this view and convince also those
Readers who might consider our reasoning too abstract and only formal
in nature, we proceed to provide some examples.

\section{Examples}

In this section we consider three examples each of them pointing out a
peculiar feature of quantum Clifford algebras and $\openZ_n$-gradings. Two
of these examples have been found by using CLIFFORD, a Maple V Rel. 5
package for quantum Clifford algebras \cite{Ablamowicz,CLIFFORD}. While
the second example is generic, the third one was taken from \cite{Fauser-vacua}
and provides an example of a physical theory which benefits
extraordinarily from using quantum Clifford algebras.

\subsection{Example 1}

This example shows that even in classical Clifford algebras one does not have
a unique access to the {\it objects\/} of the graded space. Consider
the well-known Dirac $\gamma$ matrices which generate the Dirac-Clifford
algebra $\cl_{1,3}$ and satisfy $\gamma_i\gamma_j+\gamma_j\gamma_i= 2\eta_{ij}\Id$ with
the Minkowski metric $\eta_{ij}=\diag(1,-1,-1,-1).$ The linear span of the
$\gamma$-matrices (generators) contains $1$-vectors $\bx=\sum x^i
\gamma_i.$ Define $\gamma_5=\gamma_0\gamma_1\gamma_2\gamma_3$ and note
that $\gamma_5^2=-\Id.$ If we define {\it new generators\/}
$\alpha_i:=\gamma_i\gamma_5$ which are 3-vectors(!), it is easily checked
that they nevertheless fulfill $\alpha_i\alpha_j+\alpha_j\alpha_i = 
2\eta_{ij}\Id.$ They might be called {\it vectors\/} on an equal right.

Define the map $\underline{\gamma}_5 : \cl_{1,3} \mapsto \cl_{1,3},$ 
$\bx \mapsto \bx^\prime:= \bx\gamma_5,$ lifted to $\cl_{1,3}.$ We
have thus defined two {\it different\/} Clifford maps
$\gamma : V \mapsto \cl_{1,3}$ and $\gamma^\prime : V \mapsto \cl_{1,3}$
with $\gamma^\prime := \underline{\gamma}_5 \circ \gamma.$ That is one 
can't know for sure which elements are 'vectors' even in this case.

We emphasized earlier that we did not expect the interpretation and
the mathematical aspects of classical Clifford algebras to change in such a
transformation. However, see \cite{Conradt} for a far more elaborate
application of a similar situation where {\it both\/} gradings are used. 

\subsection{Example 2}

In this example we examine the split case $\cl_{2,2} \simeq \cl_{1,1}
\otimes \cl_{1,1}$ and show the existence and irreducibility of an
$8$-dimensional representation not known in the classical representation 
theory of Clifford algebras.

We start with $\cl_{1,1}(B)$ where $B$ is given as
\begin{equation}
B := \left(\begin{array}{cc}
1 & a \\ 0 & -1 
\end{array}\right).
\end{equation}
If $a$ is zero, we have two choices for an idempotent element generating
a spinor space:
\begin{eqnarray}
\bff^-_{11}:=\frac12(\Id+\be_1), &&
\bff^+_{11}:=\frac12(\Id+\be_1\w \be_2) \, .
\end{eqnarray}
A spinor basis can be found in both cases by left multiplying by $\be_2$
which yields $\cS^\pm = <\bff^\pm_{11},\be_2 \bff^\pm_{11}>.$ The spinor spaces
$\cS^\pm$ are $2$-dimensional and the Clifford elements are represented as $2\times 2$
matrices. If $a$ is not zero, an analogous construction runs through.

Now let us put together two such algebras, as shown in \cite{Maks},
generated by $\cl_{1,1}=<\be_1,\be_2>$ and $\cl_{1,1}=<\be_3,\be_4>.$ The
bilinear form $B$ which reduces in both cases to the above setting {\it
and\/} which contains connecting elements is
\begin{equation}
B := \left(\begin{array}{cccc}
1 & a & n_{11} & n_{12} \\
0 & -1 & n_{21} & n_{22} \\
0 & 0 & 1 & a \\
0 & 0 & 0 & -1 
\end{array}\right) \, .
\end{equation}
We expect the $n_{ij}$ parameters to govern the deformation of the tensor
product in the decomposition theorem.

Searching with CLIFFORD for idempotents in this general case yields the
following fact. Let $\lambda$ be a fixed parameter. Among six choices for an 
idempotent $\bff,$ we found 
$$
\bff:=\frac12(\Id+X_1)=
\frac14(2+\lambda\,a)\Id+ \frac14\sqrt{4-\lambda^2a^2-4\lambda^2}\,{\be_1}+\frac12
\lambda \, \be_1 \wedge \be_2\,
$$
where $X_1$ is one of six different, non-trivial, and general elements $X$ in $\cl(B,V)$  satisfying $X^2 = \Id.$
This is an {\it indecomposable idempotent\/} which therefore
generates an {\it irreducible $8$ dimensional representation\/} since the
regular representation of $\cl(B,V)$ is of dimension $16.$ This fact depends on the
appearance of the non-zero $n_{ij}$ parameters. It was proved by brute force
that none of the remaining five non-trivial elements $X_i,i=2,\ldots,6,$ and
squaring to $\Id$ commuted with $X_1.$ Thus, the search showed that there is no second 
Clifford element $X_2 \neq X_1$ which would square to $\Id$ and which would commute with $X_1.$ Such an element would be necessary to decompose $\bff$ into a product 
$\bff=\prod_i \frac12 (\Id+X_i)$ where $X_iX_j=X_jX_i$ and $X_i^2=\Id.$ Since this type of
reasoning can be used to classify Clifford algebras \cite{Dimakis} we have
found a way to classify quantum Clifford algebras.

This type of an indecomposable exotic representation will occur in the
next example of a physical model and is thereby not academic.

\subsection{Example 3}

\subsubsection{Index doubling}

For a simple treatment with a computer algebra, using CLIFFORD package, and
for physical reasons not discussed here, see
\cite{Fauser-thesis,Fauser-vacua,Fauser-transition}, we introduce an index
doubling which provides us with a possibility to map the contraction and the
wedge onto a new Clifford product in the larger algebra. The benefits of
such a treatment are: the associativity of the mapped products, only one
algebra product needed during calculations, etc. 

Define the self-dual (reflexive) space $\bV = V \oplus V^*$ and
introduce generators $\be_i$ which span $V$ and $V^*$
\begin{eqnarray}
V\,=\, <\be_1,\ldots, \be_n>, &&
V^*\,=\, <\be_{n+1},\ldots, \be_{2n}> \, .
\end{eqnarray}
In this transition we require that the elements $\be_i$ from $V$ generate
a Grassmann sub-algebra and the $\be_{n+1},\ldots,\be_{2n}\in V^*$ are duals which act via
the contraction on $V.$ This gives the following conditions on the form $\bB : \bV \times \bV \mapsto \bR:$
\begin{eqnarray}
i)   && \be_i^2 = \be_i \w \be_i \w \cdot =0 \nn
ii)  && \be_{n+i}^2 = \be_{n+1} \JJbfBB \be_{n+i} \JJbfBB \cdot 
        = (\be_{n+i} \w \be_{n+i})\JJbfBB \cdot = 0 \, . 
\end{eqnarray}
Thus, with respect to the basis $<\be_1,\ldots,\be_n,\be_{n+1},\ldots,\be_{2n}>,$ $\bB$ has the following matrix:
\begin{equation}
\bB := \left(\begin{array}{cc}
0 & g \\ g^T & 0
\end{array}\right) + A \,=\, g+A,
\end{equation}
where, with an abuse of notation, the symmetric part of $\bB$ is again
denoted by $g.$ Note that we have introduced here a further freedom since
$A$ may be non-trivial also in the $V$-$V$ and $V^*$-$V^*$ sectors. This fact has certain
physical consequences which were discussed in \cite{Fauser-vacua}. The $\be_i$'s
from $V$ can be identified with Schwinger sources of quantum field theory
\cite{Fauser-transition,Fauser-thesis}.

\subsubsection{The $U(2)$-model}

We simply report here the result from \cite{Fauser-vacua} and strongly encourage
the reader to consult this work since we quote here only a
part of that work which shows the indecomposability of quantum Clifford
algebra representations and the therefrom following physical consequences.

Define $\cl(B,V) \simeq \cl_{2,2}(B)$ by specifying $\bV = <\be_i> =
<\ba_1^\dagger,\ba_2^\dagger,\ba_3,\ba_4>$ and
\begin{equation}
\bB := \frac{1}{2}\left(\begin{array}{cc}
0 & \Idm \\ \Idm & 0 
\end{array}\right) + A,
\end{equation} 
where $\Idm$ is the $2\times 2$ unit matrix and $A$ is an arbitrary but
fixed $4\times 4$ antisymmetric matrix with respect to the $\be_i$ or $\ba_i$
basis. Note furthermore that the $\ba_i$ and $\ba_i^\dagger$
fulfill the canonical anti-commutation relations, CAR, of a quantum
system: $\{ \ba_i, \ba_j^\dagger \}_+ = \delta_{ij}.$ Define
furthermore Clifford elements $N,S_i \in \openR \oplus \bV \w {\bf 
V},$ $i\in \{1,2,3\}$ such that the following relations hold:
\begin{eqnarray}
& [N,a_i]_- \,=\, -a_i,\quad [N,a^\dagger_i]_- \,=\, +a^\dagger_i
,\quad N^\dagger\,=\, N, & \nn
& [S_k,a_i]_- \,=\, \sigma_{ij}a_j,\quad \mbox{h.c.,~~~   
$k\in\{1,2,3\},$}& \nn
& [S_k,N ]_- \,= 0,\quad [S_k,S_l]_- \,=\, i\epsilon_{klm}\,S_m,
\quad S_K^\dagger = S_k \, ,&
\end{eqnarray}
where $\dagger$ is the anti-involutive map (includes a product reversion)
interchanging $\ba_i \leftrightarrow \ba_i^\dagger.$ This is the
$U(2)$ algebra if $A\equiv 0.$

Define a `vacuum', for a discussion see \cite{Fauser-vacua}, simply be
defining the expectation function --linear functional-- as the projector
onto the scalar part $<.>^A_0$ which depends now explicitly on $A.$ In
a physicist's notation $<0\mid \hat{\cH} \mid 0> \simeq $ $<\cH>^A_0$ for any
operator $\hat{\cH}$ resp. Clifford element $\cH.$

An algebraic analysis which coincides in the positive definite case with
$C^*$-algebraic results shows that this linear functional called `vacuum'
can be uniquely decomposed in certain extremal, that is indecomposable,
states. Denoting these states as spinor like $\cS_1,\cS_2$ and exotic
$\cE$ we obtain the following identity:
\begin{equation}
<.>^A_0 \;=\; \lambda_1 <.>^{\cS_1} + \lambda_2 <.>^{\cS_2}
         + \lambda_3 <.>^{\cE}, \quad
\sum \lambda_i =1.
\end{equation}
Since the regular representation of $\cl_{2,2}(B)$ is $16$ dimensional and
we find dim $\cS_1$ = dim $\cS_2 = 4,$ dim $\cE = 8$ this is a direct sum
decomposition into irreducible representations. The 'classical' case would
have led to four representations of the spinor type each $4$ dimensional.
The indecomposable exotic representation obtained from $<.>^{\cE}$ is
therefore new and it is a direct outcome of the structure of the quantum
Clifford algebra, see previous example. This representation decomposes
into two spinor like parts if $A$ vanishes identically $A\equiv 0.$

\begin{figure}[thb]
\centerline{\psfig{figure=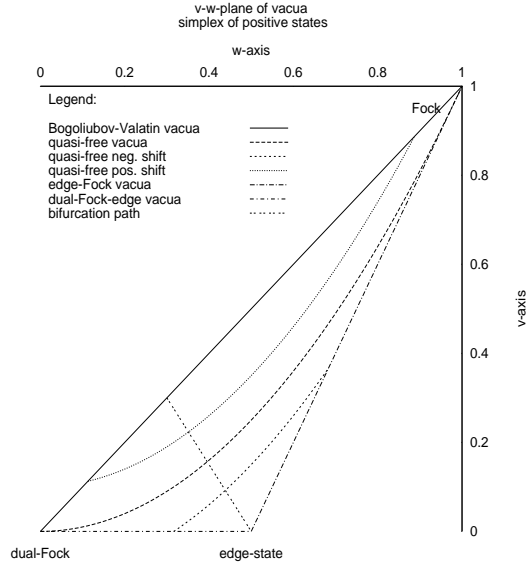,height=0.8\textwidth,width=0.53\textwidth}}
\caption{$\cl(B)$-deformation of $U(2)$ algebra}
\end{figure}

In \cite{Fauser-vacua} we obtained a $v$-$w$-plane of vacua while
implementing the $\sum \lambda_i=1$ condition and renaming of variables
into $v,w.$ There it was shown, see Figure 1, that we find free systems
of Fock and dual-Fock type which constitutes the spinor representations
$\cS_1,\,\cS_2$ and that the line connecting them contains
Bogoliubov-transformed ground-states of BCS-superconductivity. Quasi free,
that is correlation free, states are on the displayed parabola. In the
exotic state one finds spin $1$ and spin $0$ components which are beyond
Bogoliubov transformations. Every choice of $A$ fixes {\it exactly\/} one
particular state in the $v$-$w$-plane. Hence, we have solved the problem
of finding an algebraic condition on which side of the Clebsch-Gordan
identity ${\bf \frac12} \otimes {\bf \frac12} = {\bf 0} \oplus {\bf 1}$ our
algebraic system {\it has to be\/} treated.

Our model, even if only marginally discussed, shows all features we want
to see in the composite and multi-particle theory. Moreover, exotic
representations which describe 'bound objects' not capable of a
decomposition are {\it beyond\/} the treatment in \cite{Doran-states}
which mimics in Clifford algebraic terms the usual tensor method which
generically bears this problem. In this context we refer to the
interesting
work of Daviau \cite{Daviau} on de Broglie's spin fusion theory
\cite{deBroglie} and to the joint works with Stumpf and Dehnen
\cite{Fauser-iso,Fauser-positron} which are connected with algebraic
composite theories.

\section*{Acknowledgment}

The first author (BF) acknowledges a travel grant of the DFG and a critical
reading of the manuscript by Th. Konrad.

\end{document}

%% file: Makro98.tex
\newcommand{\JJ}{\mathbin{\raisebox{0.25ex}{$\footnotesize
                       \rm\vphantom{I}%
                       \_\hskip -0.25em\_%
                       \vrule width 0.6pt$}}}           

\newcommand{\w}{\wedge}
\newcommand{\bigw}{\bigwedge}
\newcommand{\dw}{\mathbin{\dot\wedge}}

\newcommand{\cE}{{\cal E}}
\newcommand{\cH}{{\cal H}}\newcommand{\cI}{{\cal I}}

\newcommand{\cS}{{\cal S}}

  \newcommand{\ba}{{\bf a}}
\newcommand{\bB}{{\bf B}}  \newcommand{\bb}{{\bf b}}

  \newcommand{\be}{{\bf e}}

\newcommand{\bR}{{\bf R}}

\newcommand{\bV}{{\bf V}}  \newcommand{\bv}{{\bf v}}
  
  \newcommand{\bx}{{\bf x}}
  \newcommand{\by}{{\bf y}}

\newcommand{\cl}{C \kern -0.1em \ell}     

\newcommand{\ut}[1]{{\setbox0=\hbox{$#1$}\mathsurround=0pt
       \rlap{\raisebox{-0.8\dp0}{\raisebox{-0.8ex}
       {\kern -0.15ex\hbox{$\tiny\sim$}\kern 0.15ex}}}#1}}
\newcommand{\uti}[1]{{\setbox0=\hbox{$#1$}\mathsurround=0pt
       \rlap{\raisebox{-0.8\dp0}{\raisebox{-0.8ex}
       {\kern -0.3ex\hbox{$\tiny\sim$}\kern 0.3ex}}}#1}}

\newdimen\arrayruleHwidth                 
     \makeatletter
     \def\Hline{\noalign{\ifnum0=`}\fi\hrule \@height \arrayruleHwidth
         \futurelet \@tempa\@xhline}
     \makeatother